\documentclass{article}
\usepackage{graphicx}

\begin{document}

\title{A statistical test for correspondence of texts to the Zipf---Mandelbrot law
\footnote{The research was supported by RFBR grant 19-51-53010.}
}
\author{Anik Chakrabarty\thanks{E-mail: a.chakrabarty@g.nsu.ru,  
Novosibirsk State University, Novosibirsk, Russia},
Mikhail Chebunin\thanks{E-mail: chebuninmikhail@gmail.com, Sobolev Institute of Mathematics, 
Novosibirsk State University, Novosibirsk, Russia},
Artyom Kovalevskii \thanks{E-mail: artyom.kovalevskii@gmail.com, Novosibirsk State Technical University, 
Novosibirsk State University, Novosibirsk, Russia},  \\
Ilya Pupyshev \thanks{E-mail: iluxa1@ngs.ru, Novosibirsk State Technical University, 
Novosibirsk State University, Novosibirsk, Russia}, 
Natalia Zakrevskaya \thanks{E-mail: natali.erlagol@gmail.com, Novosibirsk State Technical University, Novosibirsk, Russia},
Qianqian Zhou \thanks{E-mail: qianqzhou@yeah.net, School of Mathematical Sciences, Nankai University, Tianjin 300071, China}.
}
\date{}
\maketitle

\begin{abstract}

We analyse correspondence of a text to a simple probabilistic model. The model assumes that the words are selected
independently from an infinite dictionary. The probability distribution correspond to the Zipf---Mandelbrot law.
We  count sequentially the numbers of different words in the text and get the process of the numbers of different words. 
Then we estimate Zipf---Mandelbrot law parameters using the same sequence
and construct an estimate of the expectation of the number of different words in the text.
Then we subtract the corresponding values of the estimate from the sequence and normalize along the coordinate axes, obtaining
a random process on a segment from 0 to 1. We prove that this process (the empirical text bridge)
converges weakly in the uniform metric on $C (0,1)$ to a
centered Gaussian process with continuous a.s. paths. We develop and implement an algorithm
for approximate calculation of
eigenvalues of the covariance function of the limit Gaussian process, 
and then an algorithm for calculating the probability distribution
of the integral of the square of this process. We use the algorithm to analyze uniformity
of texts in English, French, Russian and Chinese.

\end{abstract}

\section{Introduction}

Our analysis is based on the fact that a text in any natural language can be divided into words.
The source material for our analysis is a text in which words are separated and punctuation marks are excluded.
In addition, all capital letters (if any) are replaced by lowercase.

We test the hypothesis $H$ that a text matches a simple probabilistic model. This model is:

1) the potential number of different words (the volume of the dictionary) is infinitely large;

2) words are selected from the dictionary independently of each other
with the same discrete probability distribution;

3) this probability distribution has a very specific form: if we arrange an infinite dictionary
according to word probabilities, then the probability of occurrence of a word with number $i$ is 
\begin{equation}\label{zm}
p_i= c (i+q)^{-\alpha}, \ \ i \ge 1, \ \ \alpha=1/\theta, \ \ 0<\theta<1, \ \  q>-1.
\end{equation} 

$c>0$ is a normalizing constant ensuring equality of the sum of the probabilities to 1.

These assumptions have a long history. Power-law character of decay of probabilities together with the infinity of the dictionary
proposed by Zipf (1936). Mandelbrot (1965) noted that  shift $q$ is required to better match real texts.
Modern large-scale studies (see Petersen et al., 2012) show that the process of the emergence of new words never stops, but
very long sequences of texts show a slightly lower frequency of new words
than this is predicted by the formula (\ref{zm}). Therefore, we will test the hypothesis $H$ for not very long texts.
(less than $ 10^5 $ words). We will demonstrate an example of non-correspondence of a long text to the hypothesis $H$.

The second assumption of the independence of the choice of consecutive words is obviously false for any meaningful text.
We can easily reject it statistically. If we calculate the relative frequency of the word {\it and} in an English text,
and then the relative frequency of the sequence {\it and and} (that is, two {\it and} in succession),
then, according to $H$, the second should be approximately equal to the square of the first.
In practice, the second is much smaller, often it is simply zero.

But we are not ready to abandon the assumption of independence. This assumption is a base of our analysis. 
Therefore, we choose a characteristic that changes little when rearranging neighboring words. 
This is the number of different words of the text, more precisely,
a sequence of numbers of different words of the text $R_1, \ldots, R_n$. Under this assumption of independence,
Bahadur proved the law of large numbers $R_n / {\bf E} R_n \to 1$ in probability, and Karlin proved
strong law of large numbers $R_n / {\bf E} R_n \to 1$ a.s. and the central limit theorem:
$(R_n - {\bf E} R_n) / \sqrt {{\bf Var} R_n}$ converges weakly to the standard normal law.

Thus, we consider a text as a random sequence, we construct the sequence 
$R_1, \ldots, R_n$ and study it using the methods of the theory of random processes: we invent parameter estimates
and construct the empirical text bridge, that is, a random process built on the parameter estimates and 
the sequence of numbers of different words.
We find the limit Gaussian process in the sense of weak convergence in $C (0,1)$. Then we calculate the distribution function $G$
of the integral of the squared limit process
based on the Smirnov (1937) formula.

Eigenvalues of the covariance function that are necessary for applying the Smirnov formula 
are calculated approximately as the eigenvalues of the matrix $Q= (q_{ij})_{i,j=1}^{100}$ composed of 
coefficients
\begin{equation}\label{qij}
q_{ij}=\int_0^1 \int_0^1 \widehat{K}(s,t) \sin \pi i s \sin \pi j t \, ds dt.
\end{equation}

We calculate $q_{ij}$ by a fast algorithm that reduces double integrals to definite integrals.

Asymptotics of $R_n$ and similar statistics (in particular, the number of unique words, 
that is, words with exactly one occurrence in the text) have been 
studied by a number of authors. The Gaussian approximation under assumptions (\ref{zm}) was studied by Karlin, 
and beyond these assumptions
by Dutko (1989), Gnedin, Hansen and Pitman (2007), Hwang and Janson (2008), Barbour and Gnedin (2009). 
Barbour (2009) proposed translated Poisson
approximation for the number of unique words. New papers by Ben-Hamou, Boucheron and Ohannessian (2017) and
Decrouez, Grabchak and Paris (2018) proved new general facts about these statistics.

The main result on which our study is based is the functional central limit theorem for the sequence
$R_1, \ldots, R_n$. It was proven by
Chebunin and Kovalevskii (2016) in preparation for this study. Note
that the Gaussian process which is the limit for this sequence
can be found also in Durieu and Wang (2016) as a limit for another prelimit process. 
Its generalization is in Durieu, Samorodnitsky and Wang (2019).

Zipf parameter estimates were proposed by Nicholls (1987), Chebunin and Kovalevskii (2019b), but we need 
a special estimate for which we can calculate
joint limit distribution of it and of the sequence $R_1, \ldots, R_n$. This estimate is proposed in
Chebunin and Kovalevskii (2019a). Zakrevskaya and Kovalevskii (2019) used the estimate in analysis of
Shakespeare's sonnets. 

Bahadur proved that under $H$ the mathematical expectation of the number of different words  
grows according to an asymptotically power law. 
This fact is known to experts in natural language processing as  Herdan's law (Herdan, 1960) or Heaps' law
(Heaps, 1978). Van Leijenhorst and van der Weide (2005), Eliazar (2011) analyzed the relationship between  
the Zipf's law and Heaps' law based
on probabilistic models that were different from $H$.

Gerlach and Altmann (2013)  noted specifically that there is no mathematically correct statistical test for correspondence 
of a text to the Zipf's law. We proposed such a test in Chebunin and Kovalevskii (2019a). We are developing
the algorithm and applying it to texts in different languages in the present paper.

The rest of the paper is organised as follows. The neccessary theoretical results are in  Section 2, 
the algorythm is in Section 3, and 
results of text analysis are in Section 4. Proofs and calculation of matrix $Q$ are in Appendixes 1 and 2.

\section{Theoretical results }

Note that constant $c$ in (\ref{zm}) is
\[
c=\left(\zeta(\alpha,q+1)\right)^{-1},
\]
\[
\zeta(\alpha,x)=\sum_{i=0}^{\infty} (i+x)^{-\alpha}
\]
is Hurvitz zeta function. 

Let $n$ be a number of words in a text.

Let  $R_k$ be the number of different words among first $k$ words of the text.

Let $R_0=0$. We have $R_1=1$, $R_0<R_1\le R_2 \le \ldots \le R_n$.

Bahadur (1960) proved that
\begin{equation}\label{Bah}
{\bf E}R_{j} \sim c^{\theta} \Gamma(1-\theta) j^{\theta},
\end{equation}
where $\Gamma(x)=\int_0^{\infty} y^{x-1} e^{-y}\,  dy$ is the Euler gamma function.
Bahadur also proved
convergence in probability
$R_j/{\bf E}R_{j} \stackrel{p}{\to} 1$.

Karlin (1967)  proved that $R_j/{\bf E} R_{j} \stackrel{a.s.}{\to} 1$, which is equivalent to 
\begin{equation}\label{Heaps} 
R_{j} \sim C_1 j^{\theta} \ {\rm a.s.,}
\end{equation}
thanks to (\ref{Bah}).
Here $C_1=c^{\theta} \Gamma(1-\theta)$, coefficients $c$ and $\theta$ are from (\ref{Bah}).

Chebunin and Kovalevskii (2016) proved the Functional Central Limit Theorem, that is, the weak convergence of the process
$\{(R_{[nt]}-{\bf E} R_{[nt]})/\sqrt{{\bf E} R_n}, \ 0\le t \le 1\}$
to a centered Gaussian
process $ Z $ with continuous a.s. sample paths and covariance function of the form
\[
K(s,t)=(s+t)^{\theta}-\max(s^{\theta}, t^{\theta}).
\]

We propose the following estimator for parameter $\theta$ (Chebunin and Kovalevskii, 2019a)
\[
\widehat{\theta}=\int_0^1  \log^+ R_{[nt]} \, d A(t)
\]
with function $A(\cdot)$ such that
\begin{equation}\label{intalog}
\int_0^1 \log t \, dA(t) =1, \ \ A(0)=A(+0)=A(1) =0,
\end{equation}
here $\log^+ x=\max(\log x, \, 0)$.
We assume $A(\cdot)$ to be the sum of a step function and a piecewise continuosly differentiable function on $[0,\, 1]$.
Let  $A(t)=0$, $t \in [0,\, \delta]$ for some $\delta\in(0,\, 1)$.

The next theorem follows from Theorems 2.1, 2.2 by Chebunin and Kovalevskii (2019a). 

{\bf Theorem 2.1} 
The estimator $\widehat{\theta}$ is strongly consistent, and
\[
\sqrt{{\bf E}R_n} (\widehat{\theta}-\theta) - \int_0^1  t^{-\theta} Z_n(t) \, dA(t) \to_p 0.
\]

From Theorem 2.1, it follows that  $\widehat{\theta}$ converges to $\theta$ with rate $({\bf E}R_n)^{-1/2}$, and 
normal random variable $ \int_0^1  t^{-\theta} Z_{\theta}(t) \, dA(t)$ has variance
$ \int_0^1 \int_0^1 (st)^{-\theta} K(s,t) \, dA(s) \, dA(t)$.

{\bf Example 2.1} Take
\[
A(t) = \left\{
\begin{array}{ll}
0, & 0 \le t \le 1/2; \\
-(\log 2)^{-1}, & 1/2<t<1; \\
0, & t=1.
\end{array}
\right.
\]

Then
\[
\widehat{\theta} = \log_2 (R_n/R_{[n/2]}), \ \ n  \ge 2.
\]

Note that, in this example, for any function $g$ on $[0,1]$, 
\[
\int_0^1 g(t)\, dA(t)=\frac{g(1)-g(1/2)}{\log 2}.
\] 

Let us introduce the process $\widehat{Z}_n$:
\[
\widehat{Z}_n(k/n)=\left(R_k-(k/n)^{\widehat{\theta}}R_n\right)/\sqrt{R_n},
\]
$0\leq k \leq n$.
Let
for $0\leq t \leq 1/n$ and $0\leq k \leq n-1$
\[
\widehat{Z}_n\left(\frac{k}{n}+t\right)=\widehat{Z}_n(k/n)+nt\left(\widehat{Z}_n((k+1)/n) - \widehat{Z}_n(k/n)\right).
\]

Let 
\[
K^0(s,t)={\bf E} Z_{\theta}^0(s) Z_{\theta}^0(t)=K(s,t)-s^{\theta}K(1,t)-t^{\theta}K(s,1)+s^{\theta}t^{\theta}K(1,1).
\]

{\bf Theorem 2.2} (Theorem 4.1 by Chebunin and Kovalevskii, 2019) 
$\widehat{Z}_n$ converges weakly to $\widehat{Z}_{\theta}$ as $n \to \infty$, where 
\[
\widehat{Z}_{\theta}(t) ={Z}_{\theta}^0(t) - t^{\theta}\log  t 
\int_0^1 u^{-\theta} Z_{\theta}(u)\,dA(u),
\]
$Z_{\theta}^0(t)=Z_{\theta}(t) - t^{\theta} Z_{\theta}(1)$, $0 \leq t \leq 1$.

The correlation function of $Z_{\theta}^0$ is given by
\[
K^0(s,t)={\bf E} Z_{\theta}^0(s) Z_{\theta}^0(t)=K(s,t)-s^{\theta}K(1,t)-t^{\theta}K(s,1)+s^{\theta}t^{\theta}K(1,1).
\]

{\bf Corollary 2.1}   
Let $\widehat{W}_n^2=\int\limits_0^1 \left(\widehat{Z}_n(t)\right)^2 dt$.
Then $\widehat{W}_n^2$ converges weakly to $\widehat{W}_{\theta}^2= \int\limits_0^1 \left(\widehat{Z}_{\theta}(t)\right)^2 dt$.

$\widehat{W}_n^2$ has the following representation 
\[
\widehat{W}_n^2=\frac{1}{3n}\sum_{k=1}^{n-1} \widehat{Z}_n\left(\frac{k}{n}\right)\left( 2\widehat{Z}_n\left(\frac{k}{n}\right)+
\widehat{Z}_n\left(\frac{k+1}{n}\right) \right).
\]

The p-value of the goodness-of fit test is $1-\widehat{F}_{\theta}(\widehat{W}_{n,obs}^2)$.
Here $\widehat{F}_{\theta}$ is the cumulative distribution function of $\widehat{W}_{\theta}^2$, and $\widehat{W}_{n,obs}^2$ 
is the observed value of $\widehat{W}_n^2$. 


One can estimate $F_{\theta}$ by simulations or find it explicitely using the  Smirnov's formula (Smirnov, 1937): 
if $W_{\theta}^2=\sum_{k=1}^{\infty}\frac{\eta^2_k}{\lambda_k}$,  
$\eta_1,\eta_2,\ldots$ are independent and have standard normal distribution, $0<\lambda_1<\lambda_2<\ldots$, then 
\begin{equation}\label{cdf}
F_{\theta}(x)=1+\frac{1}{\pi} \sum_{k=1}^{\infty} (-1)^k \int_{\lambda_{2k-1}}^{\lambda_{2k}}
\frac{e^{-\lambda x/2}}{\sqrt{-D(\lambda)}}
\cdot \frac{d\lambda}{\lambda}, \ x>0,
\end{equation}
\[
D(\lambda)=\prod_{k=1}^{\infty} \left( 1- \frac{\lambda}{\lambda_k}\right).
\]

The integrals in the RHS of (\ref{cdf}) must tend to 0 monotonically as $k \to \infty$, 
and $\lambda_k^{-1}$ are the eigenvalues of kernel $K^0$,
see Smirnov (1937).

Let us introduce the empirical bridge of a text $\widetilde{Z}_n$ by substituting
\[
r(k)=\sum_{i=1}^{\infty} \left(1-(1-\widehat{p}_i)^k\right)
\]
instead of $(k/n)^{\widehat{\theta}}R_n$ in $\widehat{Z}_n$, that is,
\[
\widetilde{Z}_n(k/n)=\left(R_k-r(k) \right)/\sqrt{R_n},
\]
$0\leq k \leq n$.

Here 
\[
\widehat{p}_i=c (i+\widehat{q})^{-1/\widehat{\theta}}, \ \ i\ge 1,
\]
$\widehat{q}$ is such that
$r(n)=R_n$.

{\bf Theorem 2.3} {\em If $H$ is true then there is $\widehat{q}$ such that $r(n)=R_n$ a.s.
There is convergence 
$\widehat{Z}_n - \widetilde{Z}_n \to 0$ a.s. uniformly on $q$ in any segment in $(-1, \infty)$.}

Proof of this theorem is postponed to Appendix 1.

\section{Algorithm}

Note that
\[
r(k) =
\sum_{i=1}^{\infty} \left(1- \sum_{j=0}^k C_k^j (-1)^j (\widehat{p}_i)^j\right)
\]
\[
= \sum_{j=1}^k C_k^j (-1)^{j+1} \sum_{i=1}^{\infty} \widehat{c}^j (i+\widehat{q})^{-j/\widehat{\theta}}
\]
\[
= \sum_{j=1}^k C_k^j (-1)^{j+1} \zeta(j/\widehat{\theta}, 1+\widehat{q}) \left(\zeta(1/\widehat{\theta}, 1+\widehat{q})\right)^{-j}.
\]

So we have the finite formula to calculate $r(1), \ldots, r(n)$. This formula is not good for large $k$ due to high 
complexity of calculations of binomial coefficients. So we use an appoximation by substituting an intergal instead of series 
residual
\[
r(k) \approx \sum_{i=1}^M (1 - (1-\widehat{p}_i)^k) + \int\limits_{M+0.5+\widehat{q}}^{\infty}
 (1-\exp(-k \widehat{c} y^{-\widehat{\alpha}}) dy
\]
\[
= \sum_{i=1}^M (1 - (1-\widehat{p}_i)^k)+(k\widehat{c})^{\widehat{\theta}} \int\limits_0^{k\widehat{c} N^{-\widehat{\alpha}}}
z^{-\widehat{\theta}} e^{-z} dz - N(1-\exp(-k  \widehat{c} N^{-\widehat{\alpha}})),
\]
$\widehat{\alpha}=1/\widehat{\theta}$, $N=M+0.5+\widehat{q}$.

We use $M=1000$. We calculate the integral using incomplete Gamma function.

We find $\widehat{q}$ by dichotomy method for $r(n) = R_n$
with 20 iterations on segment
$[-0.9, \, 40]$.

Elementary calculations give
\[
\widehat{K}(s,t)={\bf E}\widehat{Z}(s) \widehat{Z}(t)
\]
\[
= K^0(s,t)-t^{\theta} \log t \frac{K(s,1)-2^{\theta} K(s, 1/2)}{\log 2}
-s^{\theta} \log s \frac{K(t,1)-2^{\theta} K(t, 1/2)}{\log 2}
\]
\[
+s^{\theta}t^{\theta} (\log s +\log t) \frac{K(1,1)-2^{\theta} K(1, 1/2)}{\log 2}
\]
\[
+ s^{\theta}t^{\theta} \log s \log t \frac{K(1,1)-2^{\theta+1} K(1, 1/2)+2^{2\theta}K(1/2,1/2)}{\log^2 2}.
\]

Now we represent kernel $\widehat{K}$ by matrix $Q=(q_{ij})_{i,j=1}^{100}$ by calculation
\[
q_{ij}=\int_0^1 \int_0^1 \widehat{K}(s,t) \sin \pi i s \sin \pi j t \, ds dt.
\] 

Calculations are in Appendix.
We calculate eigenvalues $\widehat{\lambda}_i$, $1 \le i\le 100$, of matrix $Q$.
Then we use the Smirnov formula to calculate the p-value.

\section{Text analysis}

Let $W_n^2=\int\limits_0^1 \left(Z_n^0(t)\right)^2 dt$. It is equal to 
\[
W_n^2=\frac{1}{3n}\sum_{k=1}^{n-1} 
Z_n^0\left(\frac{k}{n}\right)\left( 2Z_n^0\left(\frac{k}{n}\right)+Z_n^0\left(\frac{k+1}{n}\right) \right).
\]

Then $W_n^2$ converges weakly to $W_{\theta}^2= \int\limits_0^1 \left(Z_{\theta}^0(t)\right)^2 dt$.

So the test rejects the basic hypothesis if $W_n^2 \geq C$. The p-value of the test is $1-F_{\theta}(W_{n,obs}^2)$.
Here $F_{\theta}$ is the cumulative distribution function of $W_{\theta}^2$ and $W_{n,obs}^2$ is a concrete value of $W_{n}^2$
 for observations under consideration.

There is an example of French poetry,
{\it Les Regrets} by Joachim du Bellay (1558), sonnet 1 (Fig. 1) and Shakespeare's sonnet 1 (Fig. 2).

\begin{figure}[htbp]
\begin{center}
\includegraphics[bb= 0 0 953 200]{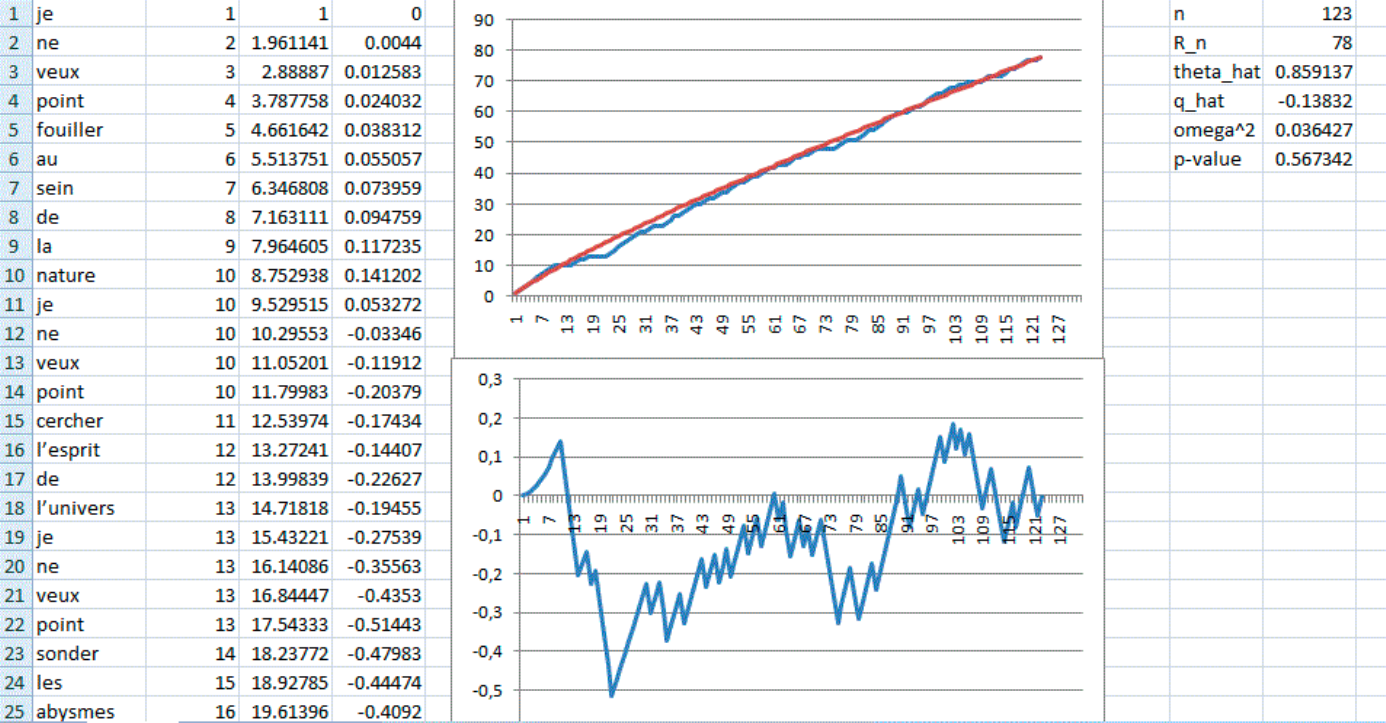}
\caption{{\it Les Regrets} by Joachim du Bellay (1558), sonnet 1}
\end{center}
\end{figure}


All results for Shakespeare's sonnets are in Appendix 3. Estimates of $\theta$ are restricted by 0.95.

\begin{figure}[htbp]
\begin{center}
\includegraphics[bb= 0 0 991 200]{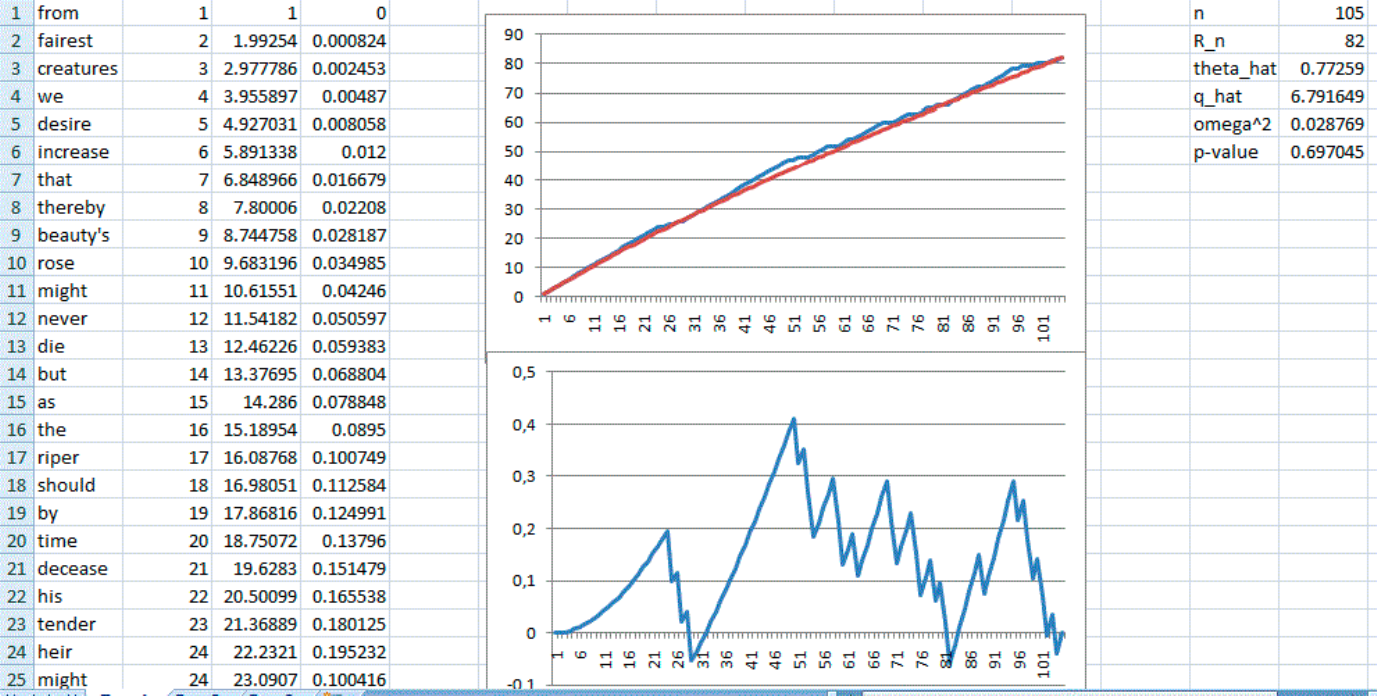}
\caption{Shakespeare's sonnet 1}
\end{center}
\end{figure}

This approach can work with texts on any language.

The first stanza of the first chapter of {\it Eugene Onegin} by Pushkin (Fig. 3).

\begin{figure}[htbp]
\begin{center}
\includegraphics[bb= 0 0 873 200]{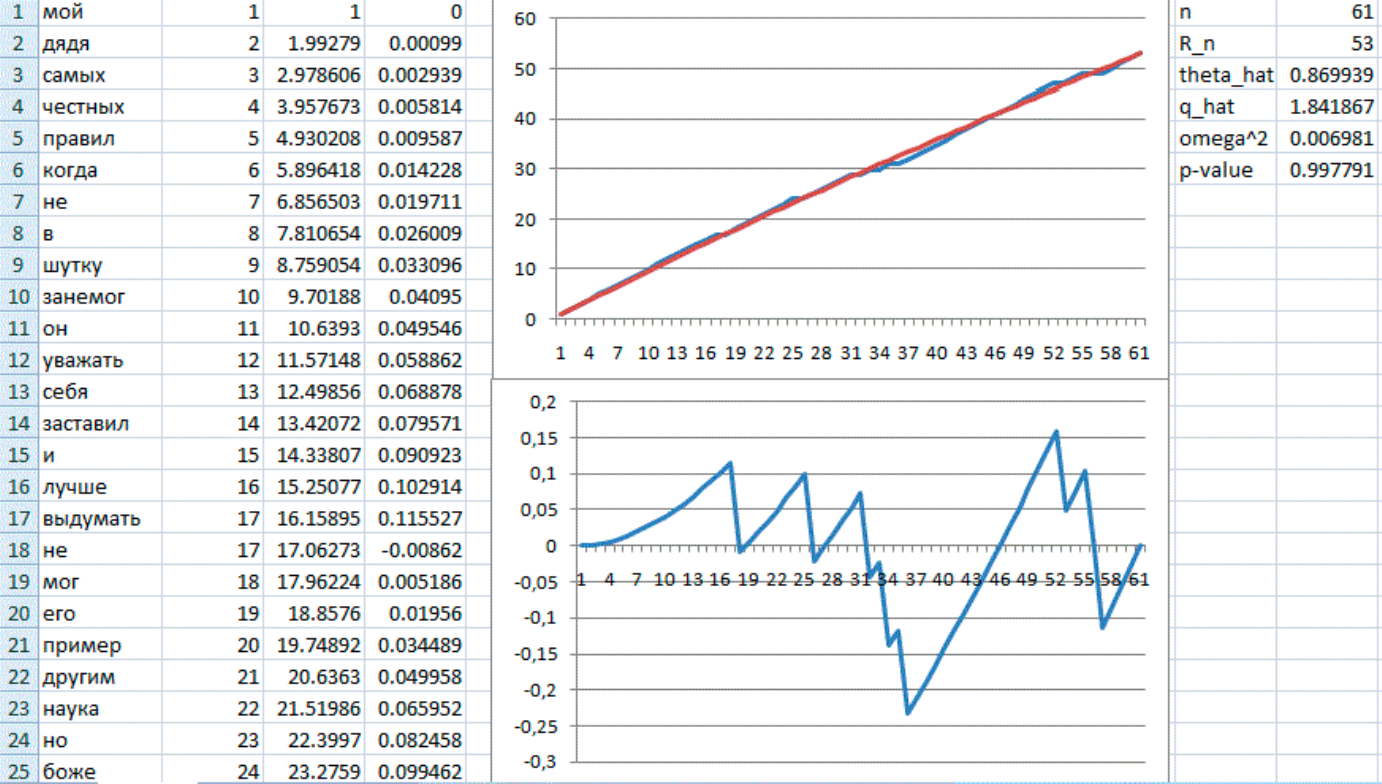}
\caption{The first stanza of the first chapter of {\it Eugene Onegin} by Pushkin}
\end{center}
\end{figure}

We used hieroglyphs instead of words when analyzing Chinese texts. 
On the first stage, we substituted hieroglyphs by their html codes.
{\it Danqing Painting} by Du Fu (Fig. 4).

\begin{figure}[htbp]
\begin{center}
\includegraphics[bb= 0 0 935 200]{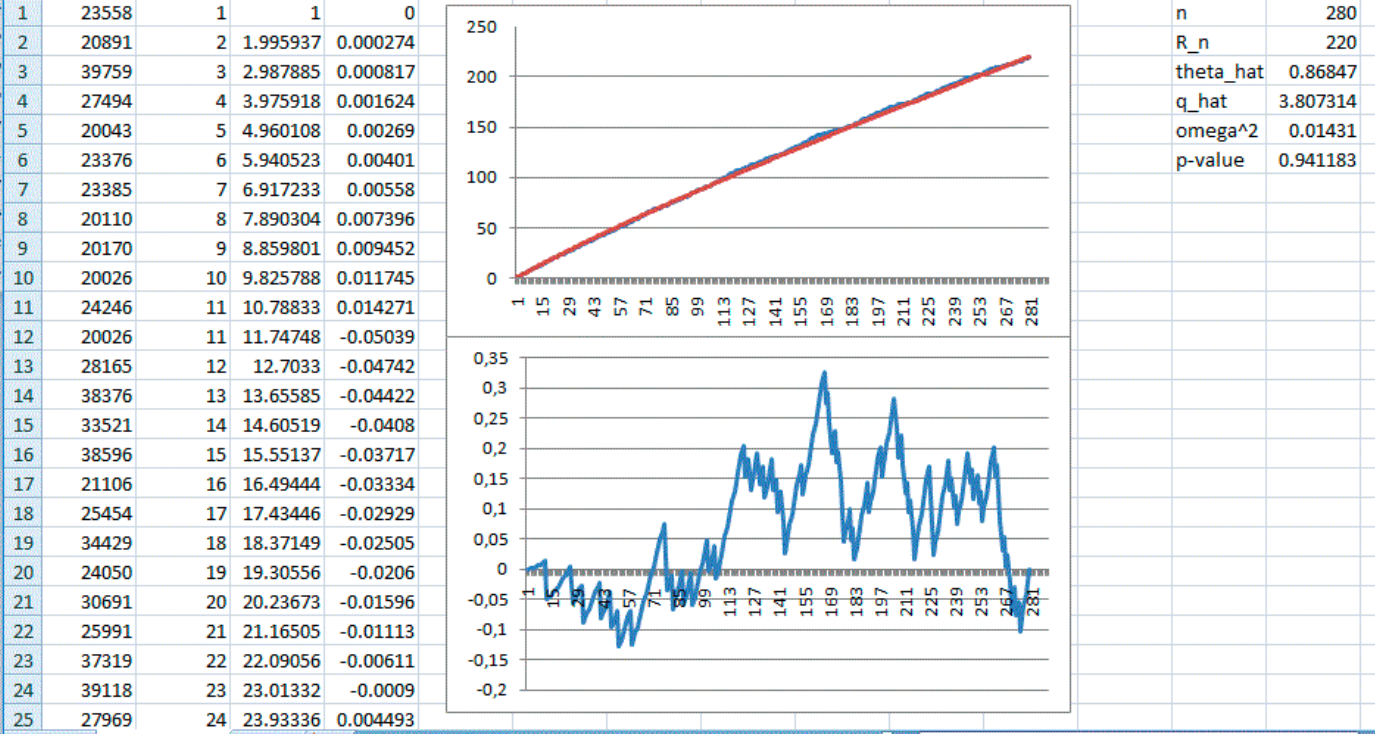}
\caption{{\it Danqing Painting} by Du Fu}
\end{center}
\end{figure}

An example of a nonhomogeneous text is Shakespeare's sonnet 1 with 3 repeated lines (Fig. 5).

\begin{figure}[htbp]
\begin{center}
\includegraphics[bb= 0 0 923 200]{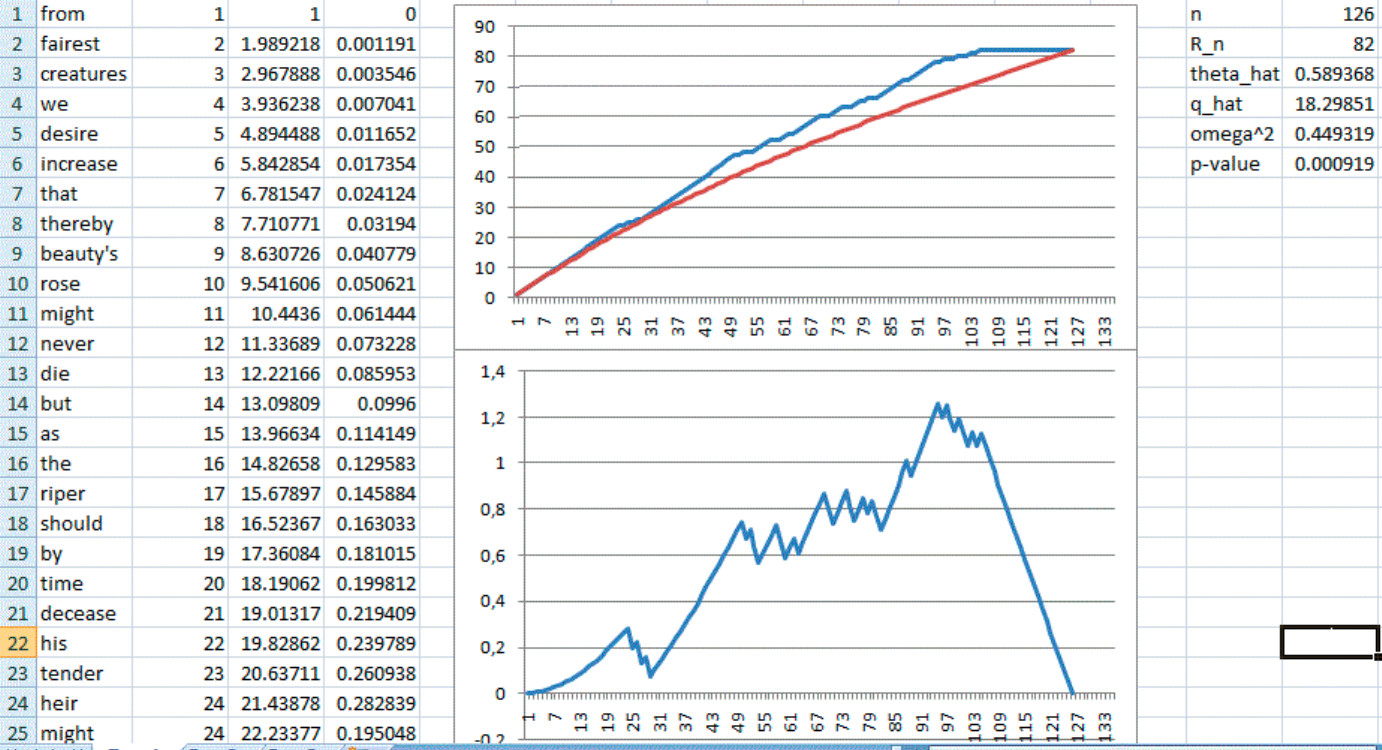}
\caption{Shakespeare's sonnet 1 with 3 repeated lines }
\end{center}
\end{figure}

Another example of nonhomogeneity is a text from different languages. 
The example is Du Bellay's sonnet 1 + Shakespeare's sonnet 1 (Fig. 6).

\begin{figure}[htbp]
\begin{center}
\includegraphics[bb= 0 0 937 200]{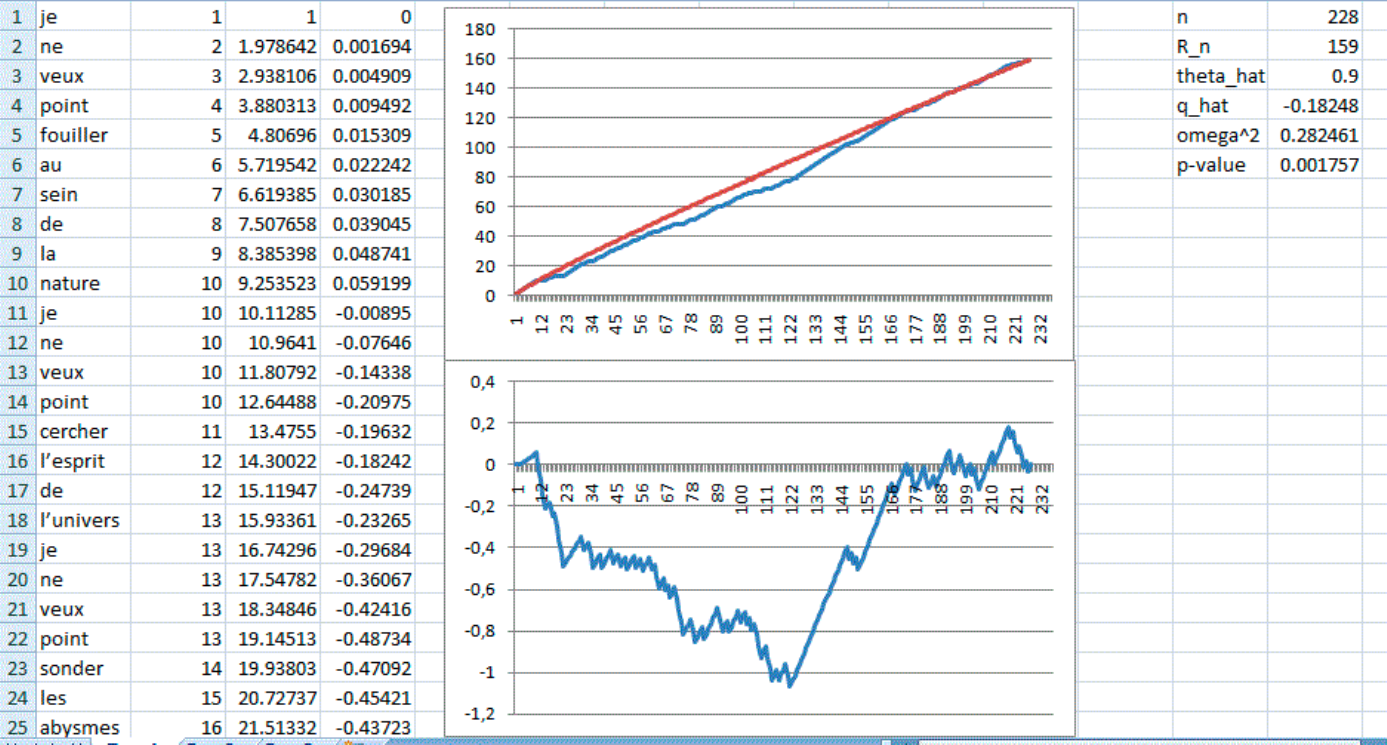}
\caption{Du Bellay's sonnet 1 + Shakespeare's sonnet 1}
\end{center}
\end{figure}

To have nonhomogeneity for a text in one language we need more long texts. Here
the first three stanzas of {\it Childe Harold's Pilgrimage}  by Byron and the first three Shakespeare's sonnets (Fig. 7).

\begin{figure}[htbp]
\begin{center}
\includegraphics[bb= 0 0 801 220]{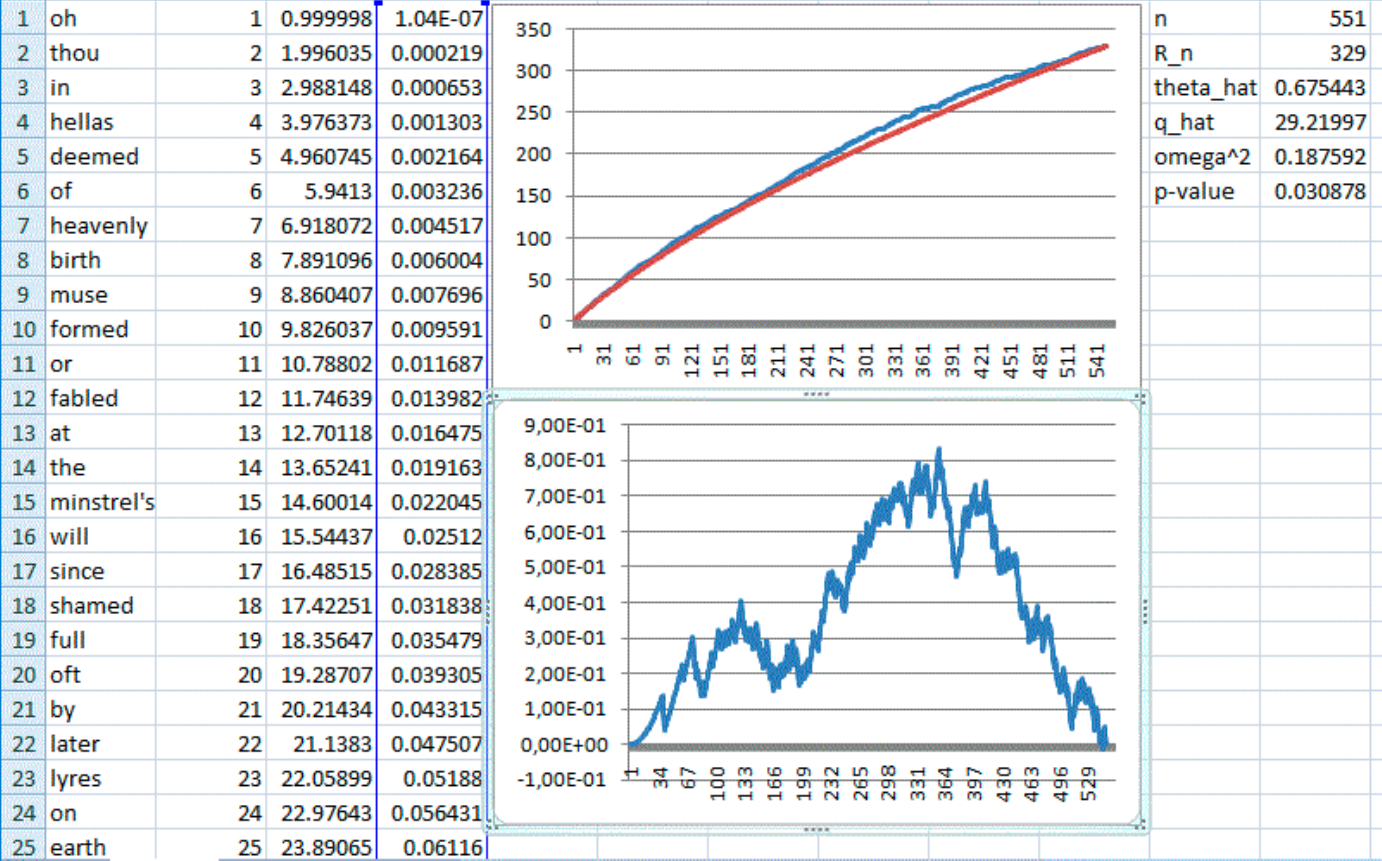}
\caption{The first three stanzas of {\it Childe Harold's Pilgrimage}  
by Byron and the first three Shakespeare's sonnets}
\end{center}
\end{figure}

\section{Discussion}

There are some open questions in the application of this approach.
The first open question is Poisson approximation. If the number of identical words is small, 
then the Gaussian approximation is inaccurate. 
Barbour (2009) proposed an approximation by a translated Poisson distribution.
We need a functional version of his theorem.

Another open question is the implementation of Simon model.
Let $R_{n,1}$ be the number of words that occur exactly once. Under the assumptions made, 
there should be convergence $R_{n,1}/R_n \to \theta$ a.s. (Karlin, 1967).
But real texts behave differently. Typically, the number of words that occur once 
is significantly less than $\widehat{\theta} R_n$.

Simon (1955)  proposed the next stochastic model: the $ (n + 1) $-th word in the text is new with  probability  $ p $; 
it coincides with each of the previous words with probability
$(1-p)/n$. 
The drawback of Simon's model is that the number of different words grows linearly.
We need some kind of hybrid of an infinite urn model and Simon's model.

\section*{Appendix 1. Proof of Theorem 2.3}

Remember that
\[
\widehat{\theta}=log_2^+ (R_n/R_{n/2}),
\]
\[ 
r(k)=\sum_{i=1}^{\infty} (1-(1-\widehat{p}_i)^k),
\]
\[
\widehat{p}_i=\widehat{c} (i + \widehat{q})^{-1/\widehat{\theta}}.
\]

From Lemma 1 in Gnedin et al. (2007) we have
\[
|{\bf E} R_n - {\bf E} R_{\Pi(n)}|<\frac{2}{n} {\bf E} R_{n,2}\le \frac{2}{n} \times \frac{n}{2}=1.
\]

Karlin (1967) proposed representation 
\[
{\bf E} R_{\Pi(t)}=\int_0^{\infty} \alpha(ty)\frac{1}{y^2} e^{-1/y} dy,
\]
\[
\alpha(x)=\max\{j|p_j\ge \frac{1}{x} \}.
\]

In our case
\[
\alpha(x)=[(cx)^{\theta}-q].
\]

We represent
\[
\Gamma(k-\gamma)=\int_0^{\infty} y^{\gamma-k-1} e^{-1/y} dy.
\]

So
\[
(ct)^{\theta} \Gamma(1-\theta)-q-1\le {\bf E} R_{\Pi(n)} \le (ct)^{\theta} \Gamma(1-\theta)-q,
\]
\[
{\bf E} R_n =(ct)^{\theta} \Gamma(1-\theta)-q+T_n, \ \ |T_n| <2,
\]
\[
|r(n)-(\widehat{c}t)^{\widehat{\theta}} \Gamma(1-\widehat{\theta})+\widehat{q}|<2 \ {\rm a.s.}
\]

Now we estimate 
\[
\max_{1\le k \le n} \frac{|(k/n)^{\widehat{\theta}} R_n - r(k) |}{\sqrt{R_n}}.
\]

We have
\[
|(k/n)^{\widehat{\theta}} R_n - r(k) |=|(k/n)^{\widehat{\theta}} r(n) - r(k) |
\]
\[
< |(k/n)^{\widehat{\theta}} ((\widehat{c}n)^{\widehat{\theta}} \Gamma(1-\widehat{\theta}) - \widehat{q} )
- (\widehat{c}k)^{\widehat{\theta}} \Gamma(1-\widehat{\theta})+ \widehat{q} |+4
\]
\[
=| \widehat{q}(1- (k/n)^{\widehat{\theta}} )|+4\le \widehat{q}+4.
\]

As $R_n \to \infty$ a.s., it is enough to prove consistency of $ \widehat{q}$ to have
\[
\max_{1\le k \le n} \frac{|(k/n)^{\widehat{\theta}} R_n - r(k) |}{\sqrt{R_n}} \to 0
\]
in probability.

\section*{Appendix 2. Calculation of $q_{ij}$}

We calculate 
\[
q_{ij}=\int_0^1 \int_0^1 \widehat{K}(s,t) \sin \pi i s \sin \pi j t \, ds dt.
\]

We know $K(s,t)=(s+t)^{\theta}-\max(s^{\theta}, t^{\theta})$,
\[
\widehat{K}(s,t)=K(s,t)-s^{\theta}K(1,t)-t^{\theta}K(s,1)+s^{\theta}t^{\theta}K(1,1)
\]
\[
-t^{\theta} \log t \frac{K(s,1)-2^{\theta} K(s, 1/2)}{\log 2}
-s^{\theta} \log s \frac{K(t,1)-2^{\theta} K(t, 1/2)}{\log 2}
\]
\[
+s^{\theta}t^{\theta} (\log s +\log t) \frac{K(1,1)-2^{\theta} K(1, 1/2)}{\log 2}
\]
\[
+ s^{\theta}t^{\theta} \log s \log t \frac{K(1,1)-2^{\theta+1} K(1, 1/2)+2^{2\theta}K(1/2,1/2)}{\log^2 2}.
\]

Let us denote
\[
J_{ij}=\int_0^1 \int_0^1 (s+t)^{\theta} \sin \pi i s \sin \pi j t \, ds dt,
\]
\[
a(\theta,k,x)=\int_0^x t^{\theta} \sin \pi k t  \, dt, \ \ b(\theta,k,x)=\int_0^x t^{\theta} \cos \pi k t  \, dt,
\]
\[
A_i = \int_0^1 t^{\theta} \sin \pi i t  \, dt = a(\theta,i,1), \ \ 
B_i = \int_1^2 t^{\theta} \sin \pi i t  \, dt = a(\theta,i,2)-a(\theta,i,1),
\]
\[
C_i = \int_0^1 t^{\theta+1} \cos \pi i t  \, dt = b(\theta+1,i,1), 
\]
\[
D_i = \int_1^2 t^{\theta}(2-t) \cos \pi i t  \, dt = 2 b(\theta,i,2) - 2 b(\theta,i,1) - b(\theta+1,i,2)+b(\theta+1,i,1),
\]
\[
E_{ij} = \int_0^1 t^{\theta} \sin \pi i t \cos \pi j t \, dt, \ \ F_i = \int_0^1 K(t,1)  \sin \pi i t  \, dt,
\]
\[
G_i = \int_0^1 t^{\theta} \log t \sin \pi i t  \, dt, \ \ H_i = \int_0^1 K(t,1/2)  \sin \pi i t  \, dt.
\]

Then we have
\[
q_{ij}=J_{ij}-\frac{1}{\pi j}(A_i - E_{ij}) -  \frac{1}{\pi i}(A_j- E_{ji}) 
\]
\[
- A_i F_j - A_j F_i + K(1,1) A_i A_j
\]
\[
-G_j \frac{F_i-2^{\theta} H_i}{\log 2}
-G_i \frac{F_j-2^{\theta} H_j}{\log 2}\]
\[
+(A_i G_j+A_j G_i) \frac{K(1,1)-2^{\theta} K(1, 1/2)}{\log 2}
\]
\[
+ G_i G_j \frac{K(1,1)-2^{\theta+1} K(1, 1/2)+2^{2\theta}K(1/2,1/2)}{\log^2 2}.
\]

We calculate $J_{ij}$ substituting $t=u-s$, $s \le u \le s+1$. So
\[
J_{ij}=
\int_0^1 u^{\theta} \, du \int_0^u  \sin \pi i s \sin \pi j (u-s) \, ds +
\int_1^2 u^{\theta} \, du \int_{u-1}^1  \sin \pi i s \sin \pi j (u-s) \, ds.
\]

If $ i \ne j$ then we have
\[
J_{ij} = \frac{1}{2} \int_0^1 u^{\theta} \, du \int_0^u  (\cos \pi( i s - ju + js) - \cos \pi( i s + ju - js)) \, ds 
\]
\[
+
\frac{1}{2} \int_1^2 u^{\theta} \, du \int_{u-1}^1  (\cos \pi( i s - ju + js) - \cos \pi( i s + ju - js)) \, ds 
\]
\[
= \frac{1}{2} \int_0^1 u^{\theta}  \left(\frac{\sin \pi i u + \sin \pi j u}{\pi(i+j)} - 
\frac{\sin \pi i u - \sin \pi j u}{\pi(i-j)}\right) \, du   
\]
\[
+
\frac{1}{2} \int_1^2  u^{\theta}  \left(\frac{\sin (\pi (i+j) -\pi j u) + \sin (\pi (i+j) -\pi i u)}{\pi(i+j)} - 
\frac{\sin (\pi (i-j) +\pi j u) + \sin (\pi (i-j) -\pi i u)}{\pi(i-j)}\right) \, du  
\]
\[
= \frac{i A_j - j A_i -(-1)^{i+j}(i B_j - j B_i)}{\pi(i^2-j^2)}.
\]

If $ i =j$ then we have
\[
J_{ii}=\frac{1}{2} \int_0^1 u^{\theta} \, du \int_0^u  (\cos \pi( 2i s - iu) - \cos \pi iu ) \, ds +
\frac{1}{2} \int_1^2 u^{\theta} \, du \int_{u-1}^1  (\cos \pi( 2i s - iu) - \cos \pi i u ) \, ds 
\]
\[
= \frac{1}{2} \int_0^1 \left( u^{\theta}  \frac{\sin \pi i u }{\pi i} - u^{\theta+1}  \cos \pi i u \right) \, du   
\]
\[
+
\frac{1}{2} \int_1^2  u^{\theta}  \frac{\sin (2\pi i -\pi i u) - \sin (2\pi iu - 2\pi i -\pi i u)}{2\pi i} \, du - 
\frac{1}{2} \int_1^2 u^{\theta}(2-u) \cos \pi i u  \, du
\]
\[
= \frac{A_i -  B_i}{2\pi i}- \frac{C_i+D_i}{2}.
\]

We know 
\[
a(\theta,k,x)=\pi k x^{\theta+2} \, {_1 F_2} \left(\theta/2+1, 3/2, \theta/2+2, -k^2 \pi^2 x^2/4 \right)/(\theta+2), 
\]
\[
b(\theta,k,x)=x^{\theta+1} \, {_1 F_2} \left(\theta/2+1/2, 1/2, \theta/2+3/2, -k^2 \pi^2 x^2/4 \right)/(\theta+1), 
\]
\[
E_{ij} = \frac{1}{2} A_{i+j}+\frac{1}{2} A_{i-j},
\]
\[
F_i=(-1)^i B_i + \frac{(-1)^i-1}{\pi i},
\]
\[
G_i = \int_0^1 t^{\theta} \log t \sum_{k=0}^{\infty} \frac{(-1)^k (\pi i t)^{2k+1}}{(2k+1)!}  \, dt
= -\sum_{k=0}^{\infty} \frac{(-1)^k (\pi i )^{2k+1}}{(2k+1)!(2k+\theta+2)^2}
\]
\[
= - {\pi i} \, {_2F_3} \left(\theta/2+1, \theta/2+1, 3/2, \theta/2+2, \theta/2+2, -i^2 \pi^2/4 \right)/(\theta+2)^2,
\] 
\[
H_i = \int_{1/2}^{3/2} t^{\theta}  \sin (\pi i t - \pi i/2) \, dt + (\cos(\pi i/2)-1)2^{-\theta}/(\pi i)
- \int_{1/2}^1 t^{\theta} \sin \pi i t  \, dt
\]
\[
= \pi i \cos (\pi i /2) \bigg( (3/2)^{\theta+2} {_1 F_2} \left(\theta/2+1, 3/2, \theta/2+2, -9i^2 \pi^2/16 \right) 
\]
\[
- 2^{-\theta-2} {_1 F_2} \left(\theta/2+1, 3/2, \theta/2+2, -i^2 \pi^2/16 \right) \bigg)/(\theta+2)
\]
\[
- \sin (\pi i /2) \bigg( (3/2)^{\theta+1} {_1 F_2} \left(\theta/2+1/2, 1/2, \theta/2+3/2, -9i^2 \pi^2/16 \right) 
\]
\[
- 2^{-\theta-1} {_1 F_2} \left(\theta/2+1/2, 1/2, \theta/2+3/2, -i^2 \pi^2/16 \right) \bigg)/(\theta+1)
\]
\[
+(\cos(\pi i/2)-1)2^{-\theta}/(\pi i) - A_i + \pi i 2^{-\theta-2} {_1 F_2} 
\left(\theta/2+1, 3/2, \theta/2+2, -i^2 \pi^2/16 \right)/(\theta+2).
\]

$_1 F_2$ and $_2 F_3$ are generalized hypergeometric functions, hyp1f2 and hyp2f3 in Python.

\section*{Appendix 3. Shakespeare's sonnets}

\begin{tabular}{|c|c|c|c|c|c|c|} \hline
           &     &       &                    &                &              &       \\
Sonnet No. & $n$ & $R_n$ & $\widehat{\theta}$ &  $\widehat{q}$ & $\omega^2$ & p-value \\ \hline
I	&	98	&	77	&	0.7911	&	5.1473	&	0.03139	&	0.6646 \\
II	&	106	&	82	&	0.7877	&	5.2669	&	0.04696	&	0.4492 \\
III	&	106	&	75	&	0.7859	&	2.9632	&	0.01184	&	0.9769 \\
IV	&	96	&	74	&	0.8344	&	2.0164	&	0.02212	&	0.8257 \\
V	&	96	&	80	&	0.9469	&	-0.5618	&	0.01114	&	0.9824 \\
VI	&	100	&	70	&	0.7030	&	7.4855	&	0.03214	&	0.6579 \\
VII	&	92	&	72	&	0.9031	&	-0.0988	&	0.00808	&	0.9970 \\
VIII	&	101	&	78	&	0.8931	&	0.0952	&	0.01886	&	0.8818 \\
IX	&	110	&	82	&	0.7877	&	4.2580	&	0.02714	&	0.7374 \\
X	&	106	&	75	&	0.8890	&	-0.2561	&	0.00732	&	0.9985 \\
XI	&	107	&	81	&	0.8804	&	0.3084	&	0.05376	&	0.3641 \\
XII	&	109	&	81	&	0.9136	&	-0.4344	&	0.03236	&	0.6379 \\
XIII	&	100	&	72	&	0.9031	&	-0.4163	&	0.01530	&	0.9375 \\
XIV	&	106	&	78	&	0.8097	&	2.5711	&	0.02065	&	0.8525 \\
XV	&	104	&	79	&	0.8280	&	2.2487	&	0.03427	&	0.6150 \\
XVI	&	101	&	80	&	0.8301	&	3.0005	&	0.02032	&	0.8580 \\
XVII	&	115	&	89	&	0.8319	&	2.7590	&	0.04821	&	0.4288 \\
XVIII	&	107	&	79	&	0.8119	&	2.5560	&	0.02585	&	0.7592 \\
XIX	&	106	&	83	&	0.9500	&	-0.7264	&	0.02625	&	0.7437 \\
XX	&	106	&	85	&	0.8858	&	0.7206	&	0.00733	&	0.9984 \\
XXI	&	109	&	78	&	0.8260	&	1.4855	&	0.01187	&	0.9768 \\
XXII	&	112	&	74	&	0.7665	&	2.8193	&	0.02753	&	0.7315 \\
XXIII	&	105	&	77	&	0.8074	&	2.5866	&	0.05756	&	0.3420 \\
XXIV	&	112	&	78	&	0.7776	&	3.2625	&	0.02204	&	0.8286 \\
XXV	&	101	&	78	&	0.8591	&	1.1303	&	0.01408	&	0.9537 \\
XXVI	&	111	&	76	&	0.8217	&	1.1439	&	0.04744	&	0.4388 \\
XXVII	&	103	&	77	&	0.8405	&	1.4188	&	0.03409	&	0.6166 \\
XXVIII	&	109	&	75	&	0.8713	&	-0.0671	&	0.01646	&	0.9210 \\
XXIX	&	107	&	79	&	0.9462	&	-0.7709	&	0.01810	&	0.8933 \\
XXX	&	109	&	81	&	0.7549	&	6.4163	&	0.02698	&	0.7415 \\
XXXI	&	108	&	73	&	0.6200	&	14.3674	&	0.14907	&	0.0698 \\
XXXII	&	104	&	80	&	0.9125	&	-0.3063	&	0.03731	&	0.5592 \\
XXXIII	&	103	&	82	&	0.8981	&	0.2218	&	0.00551	&	0.9998 \\
XXXIV	&	112	&	86	&	0.8564	&	1.3410	&	0.01778	&	0.9004 \\
XXXV	&	99	&	74	&	0.8171	&	2.3136	&	0.01207	&	0.9750 \\
XXXVI	&	107	&	76	&	0.7561	&	4.7499	&	0.07452	&	0.2364 \\
XXXVII	&	107	&	74	&	0.7500	&	4.4181	&	0.04462	&	0.4813 \\
XXXVIII	&	103	&	83	&	0.8205	&	4.2305	&	0.01956	&	0.8713 \\
XXXIX	&	111	&	77	&	0.8405	&	0.7001	&	0.02212	&	0.8255 \\
XL	&	111	&	80	&	0.9500	&	-0.8142	&	0.21140	&	0.0068 \\
\hline
\end{tabular}

\begin{tabular}{|c|c|c|c|c|c|c|} \hline
           &     &       &                    &                &              &       \\
Sonnet No. & $n$ & $R_n$ & $\widehat{\theta}$ &  $\widehat{q}$ & $\omega^2$ & p-value \\ \hline
XLI	&	102	&	72	&	0.7949	&	2.3792	&	0.06757	&	0.2699 \\
XLII	&	120	&	71	&	0.7235	&	3.2051	&	0.04893	&	0.4369 \\
XLIII	&	114	&	69	&	0.5850	&	13.1416	&	0.14373	&	0.0823 \\
XLIV	&	110	&	83	&	0.8359	&	1.9066	&	0.01353	&	0.9602 \\
XLV	&	98	&	78	&	0.9500	&	-0.7064	&	0.01251	&	0.9706 \\
XLVI	&	107	&	67	&	0.8182	&	0.4581	&	0.05119	&	0.3994 \\
XLVII	&	113	&	68	&	0.8021	&	0.6272	&	0.02254	&	0.8192 \\
XLVIII	&	111	&	83	&	0.8515	&	1.1455	&	0.01952	&	0.8713 \\
XLIX	&	104	&	74	&	0.8001	&	2.3442	&	0.01550	&	0.9353 \\
L	&	110	&	83	&	0.8359	&	1.9066	&	0.01981	&	0.8667 \\
LI	&	108	&	80	&	0.8301	&	1.8074	&	0.03863	&	0.5495 \\
LII	&	102	&	76	&	0.8050	&	2.9384	&	0.00851	&	0.9958 \\
LIII	&	99	&	62	&	0.7063	&	4.2012	&	0.03664	&	0.5900 \\
LIV	&	103	&	77	&	0.7749	&	4.8369	&	0.01686	&	0.9156 \\
LV	&	99	&	74	&	0.7832	&	4.1136	&	0.04431	&	0.4807 \\
LVI	&	105	&	80	&	0.8625	&	0.8972	&	0.01015	&	0.9890 \\
LVII	&	110	&	79	&	0.8944	&	-0.2682	&	0.01576	&	0.9311 \\
LVIII	&	103	&	70	&	0.8439	&	0.3739	&	0.01275	&	0.9685 \\
LIX	&	101	&	78	&	0.8260	&	2.6102	&	0.01904	&	0.8801 \\
LX	&	100	&	79	&	0.7960	&	5.1082	&	0.01822	&	0.8941 \\
LXI	&	114	&	78	&	0.8097	&	1.6042	&	0.01677	&	0.9168 \\
LXII	&	98	&	66	&	0.8946	&	-0.4938	&	0.01837	&	0.8899 \\
LXIII	&	103	&	78	&	0.7618	&	6.2212	&	0.02466	&	0.7823 \\
LXIV	&	104	&	74	&	0.9241	&	-0.6534	&	0.06417	&	0.2700 \\
LXV	&	106	&	85	&	0.9337	&	-0.4770	&	0.01424	&	0.9511 \\
LXVI	&	82	&	65	&	0.8524	&	1.3890	&	0.00781	&	0.9976 \\
LXVII	&	100	&	74	&	0.9241	&	-0.5877	&	0.01644	&	0.9204 \\
LXVIII	&	99	&	75	&	0.9500	&	-0.7739	&	0.02632	&	0.7423 \\
LXIX	&	114	&	81	&	0.8480	&	0.7358	&	0.01090	&	0.9843 \\
LXX	&	102	&	75	&	0.8890	&	-0.1060	&	0.01452	&	0.9481 \\
LXXI	&	115	&	76	&	0.7885	&	1.9740	&	0.01484	&	0.9442 \\
LXXII	&	110	&	75	&	0.7370	&	4.9806	&	0.03965	&	0.5455 \\
LXXIII	&	112	&	77	&	0.7590	&	3.9825	&	0.01696	&	0.9140 \\
LXXIV	&	106	&	67	&	0.6738	&	6.6998	&	0.10606	&	0.1348 \\
LXXV	&	109	&	77	&	0.8074	&	2.0411	&	0.02854	&	0.7118 \\
LXXVI	&	107	&	72	&	0.6464	&	11.2969	&	0.10129	&	0.1534 \\
LXXVII	&	101	&	70	&	0.9500	&	-0.8438	&	0.02176	&	0.8277 \\
LXXVIII	&	102	&	74	&	0.7665	&	4.4436	&	0.07051	&	0.2570 \\
LXXIX	&	108	&	74	&	0.7500	&	4.2387	&	0.03531	&	0.6059 \\
LXXX	&	108	&	81	&	0.8641	&	0.7063	&	0.01128	&	0.9815 \\
\hline
\end{tabular}

\begin{tabular}{|c|c|c|c|c|c|c|} \hline
           &     &       &                    &                &              &       \\
Sonnet No. & $n$ & $R_n$ & $\widehat{\theta}$ &  $\widehat{q}$ & $\omega^2$ & p-value \\ \hline
LXXXI	&	107	&	77	&	0.8074	&	2.2956	&	0.00920	&	0.9934 \\
LXXXII	&	99	&	79	&	0.9115	&	-0.1166	&	0.01598	&	0.9276 \\
LXXXIII	&	108	&	79	&	0.8775	&	0.1367	&	0.01685	&	0.9150 \\
LXXXIV	&	106	&	77	&	0.8074	&	2.4361	&	0.01776	&	0.9015 \\
LXXXV	&	103	&	80	&	0.7370	&	9.6313	&	0.04676	&	0.4585 \\
LXXXVI	&	101	&	80	&	0.8957	&	0.2286	&	0.01423	&	0.9517 \\
LXXXVII	&	110	&	81	&	0.9305	&	-0.6463	&	0.02146	&	0.8340 \\
LXXXVIII	&	103	&	74	&	0.7176	&	7.5928	&	0.07136	&	0.2616 \\
LXXXIX	&	104	&	79	&	0.8280	&	2.2487	&	0.03930	&	0.5403 \\
XC	&	111	&	74	&	0.7176	&	5.6453	&	0.02559	&	0.7669 \\
XCI	&	107	&	76	&	0.9500	&	-0.8264	&	0.09217	&	0.1262 \\
XCI	&	112	&	74	&	0.7337	&	4.4914	&	0.03364	&	0.6326 \\
XCIII	&	107	&	79	&	0.8443	&	1.1796	&	0.01119	&	0.9822 \\
XCIV	&	96	&	74	&	0.9500	&	-0.7554	&	0.03898	&	0.5299 \\
XCV	&	104	&	82	&	0.9500	&	-0.7170	&	0.04642	&	0.4318 \\
XCVI	&	111	&	82	&	0.8657	&	0.5435	&	0.02981	&	0.6856 \\
XCVII	&	100	&	84	&	0.9166	&	0.1683	&	0.00588	&	0.9997 \\
XCVIII	&	110	&	79	&	0.7960	&	2.8866	&	0.01799	&	0.8978 \\
XCIX	&	115	&	81	&	0.7853	&	3.1985	&	0.02053	&	0.8551 \\
C	&	101	&	78	&	0.9279	&	-0.5258	&	0.00548	&	0.9998 \\
CI	&	106	&	77	&	0.7590	&	5.1635	&	0.01166	&	0.9782 \\
CII	&	108	&	79	&	0.8443	&	1.0814	&	0.00717	&	0.9987 \\
CIII	&	108	&	74	&	0.6095	&	16.4318	&	0.21617	&	0.0247 \\
CIV	&	107	&	85	&	0.9500	&	-0.7013	&	0.01698	&	0.9115 \\
CV	&	98	&	62	&	0.6874	&	5.3704	&	0.02029	&	0.8592 \\
CVI	&	103	&	75	&	0.9069	&	-0.4228	&	0.02428	&	0.7829 \\
CVII	&	104	&	83	&	0.8672	&	1.3581	&	0.00865	&	0.9955 \\
CVIII	&	106	&	85	&	0.9016	&	0.2122	&	0.02197	&	0.8261 \\
CIX	&	108	&	72	&	0.7776	&	2.3672	&	0.01819	&	0.8947 \\
CX	&	110	&	84	&	0.8225	&	2.8207	&	0.01456	&	0.9479 \\
CXI	&	103	&	76	&	0.7885	&	3.6411	&	0.02730	&	0.7344 \\
CXII	&	109	&	69	&	0.7510	&	2.7209	&	0.03853	&	0.5594 \\ 
CXIII	&	110	&	74	&	0.7665	&	3.0736	&	0.03382	&	0.6274 \\
CXIV	&	106	&	75	&	0.6897	&	9.5564	&	0.10235	&	0.1416 \\
CXV	&	105	&	83	&	0.7026	&	15.1808	&	0.09517	&	0.1605 \\
CXVI	&	103	&	78	&	0.8260	&	2.2610	&	0.01148	&	0.9800 \\
CXVII	&	103	&	73	&	0.9044	&	-0.4676	&	0.02718	&	0.7296 \\
CXVIII	&	105	&	80	&	0.8625	&	0.8972	&	0.03203	&	0.6482 \\
CXIX	&	106	&	77	&	0.8238	&	1.7197	&	0.05184	&	0.3919 \\
CXX	&	108	&	76	&	0.8050	&	2.0566	&	0.02387	&	0.7950 \\
\hline
\end{tabular}

\begin{tabular}{|c|c|c|c|c|c|c|} \hline
           &     &       &                    &                &              &       \\
Sonnet No. & $n$ & $R_n$ & $\widehat{\theta}$ &  $\widehat{q}$ & $\omega^2$ & p-value \\ \hline
CXXI	&	107	&	79	&	0.7188	&	9.0289	&	0.03333	&	0.6383 \\
CXXII	&	100	&	80	&	0.8141	&	4.2738	&	0.02785	&	0.7235 \\
CXXIII	&	109	&	80	&	0.8301	&	1.6843	&	0.01845	&	0.8899 \\
CXXIV	&	102	&	76	&	0.8729	&	0.3330	&	0.01468	&	0.9461 \\
CXXV	&	100	&	83	&	0.8672	&	2.0726	&	0.00523	&	0.9999 \\
CXXVI	&	88	&	68	&	0.8207	&	2.4595	&	0.02146	&	0.8380 \\
CXXVII	&	102	&	79	&	0.8608	&	1.1231	&	0.02098	&	0.8454 \\
CXXVIII	&	104	&	80	&	0.8790	&	0.4621	&	0.03259	&	0.6373 \\
CXXIX	&	102	&	71	&	0.7574	&	3.9305	&	0.02730	&	0.7358 \\
CXXX	&	114	&	81	&	0.8641	&	0.2990	&	0.00732	&	0.9985 \\
CXXXI	&	111	&	73	&	0.7636	&	2.8441	&	0.03039	&	0.6829 \\
CXXXII	&	107	&	80	&	0.7984	&	3.5647	&	0.01857	&	0.8884 \\
CXXXIII	&	112	&	73	&	0.6820	&	7.4823	&	0.03948	&	0.5530 \\
CXXXIV	&	112	&	71	&	0.8099	&	0.7968	&	0.05417	&	0.3717 \\
CXXXV	&	108	&	71	&	0.7068	&	5.7796	&	0.11681	&	0.1039 \\
CXXXVI	&	114	&	68	&	0.5956	&	11.5040	&	0.13579	&	0.0915 \\
CXXXVII	&	115	&	81	&	0.8480	&	0.6669	&	0.00795	&	0.9973 \\
CXXXVIII	&	109	&	74	&	0.7832	&	2.4423	&	0.01640	&	0.9224 \\
CXXXIX	&	113	&	82	&	0.9500	&	-0.8080	&	0.02000	&	0.8600 \\
CXL	&	111	&	84	&	0.7776	&	5.4863	&	0.01730	&	0.9088 \\
CXLI	&	110	&	81	&	0.8007	&	3.1756	&	0.01889	&	0.8830 \\
CXLII	&	109	&	74	&	0.7176	&	6.0534	&	0.02093	&	0.8485 \\
CXLIII	&	112	&	78	&	0.7155	&	7.2400	&	0.04597	&	0.4701 \\
CXLIV	&	105	&	74	&	0.8171	&	1.5352	&	0.00950	&	0.9922 \\
CXLV	&	89	&	70	&	0.8439	&	1.7819	&	0.01051	&	0.9869 \\
CXLVI	&	105	&	76	&	0.8217	&	1.7322	&	0.01653	&	0.9205 \\
CXLVII	&	100	&	72	&	0.7949	&	2.6716	&	0.01535	&	0.9375 \\
CXLVIII	&	116	&	77	&	0.8745	&	-0.2326	&	0.05124	&	0.3900 \\
CXLIX	&	110	&	72	&	0.9220	&	-0.7360	&	0.02077	&	0.8471 \\
CL	&	112	&	72	&	0.7437	&	3.4057	&	0.06262	&	0.3135 \\
CLI	&	110	&	82	&	0.7877	&	4.2580	&	0.03491	&	0.6088 \\
CLII	&	115	&	70	&	0.7030	&	4.5331	&	0.03797	&	0.5716 \\ 
CLIII	&	102	&	75	&	0.8194	&	2.0017	&	0.06991	&	0.2506 \\
CLIV	&	100	&	74	&	0.7500	&	6.0394	&	0.03911	&	0.5516 \\
\hline
\end{tabular}

\bigskip

{\bf Acknowledgements}

The research was supported by RFBR grant 19-51-53010.
The authors like to thank Sergey Foss 
for helpful and constructive comments and suggestions.

\bigskip


\bigskip

\footnotesize

{\sc Bahadur, R. R.}, 1960. On the number of distinct values in a large sample from an infinite discrete distribution.
Proceedings of the National Institute of Sciences of India, 26A, Supp. II, 67--75.

{\sc Barbour, A. D.}, 2009. Univariate approximations in the infinite occupancy
scheme. Alea 6, 415--433.

{\sc Barbour, A. D.,  Gnedin, A. V.}, 2009.
Small counts in the infinite occupancy scheme.
Electronic Journal of Probability, 
Vol. 14, Paper no. 13, 365--384.

{\sc Ben-Hamou, A., Boucheron,  S., Ohannessian,  M. I.}, 2017. 
{Concentration inequalities in the infinite urn scheme for occupancy counts and the missing mass, with applications},
Bernoulli, V. 23, No. 1, 249--287.

{\sc Chebunin, M., Kovalevskii, A.}, 2016. 
Functional central limit theorems for certain statistics in an infinite urn scheme. 
Statistics and Probability Letters, V. 119, 344--348.

{\sc Chebunin, M. G., Kovalevskii, A. P.}, 2019a.
A statistical test for the Zipf's law by deviations from the Heaps' law, Siberian Electronic Mathematical Reports, 
V. 16, 1822--1832.

{\sc Chebunin, M., Kovalevskii, A.}, 2019b. 
Asymptotically Normal Estimators for Zipf's Law, 
Sankhya A, V. 81, 482--492.

{\sc Decrouez, G., Grabchak, M., Paris,   Q.},  2018.
{Finite sample properties of the mean occupancy counts and probabilities},
Bernoulli, V. 24, No. 3, 1910--1941.

{\sc Durieu, O.,  Wang, Y.}, 2016. 
{From infinite urn schemes to decompositions of self-similar Gaussian processes}, 
Electron. J. Probab., V. 21, paper No. 43, 23 pp.

{\sc Durieu, O., Samorodnitsky, G., Wang, Y.}, 2019.
{From infinite urn schemes to self-similar stable processes},
Stochastic Processes and their Applications (in press).

{\sc Dutko, M.}, 1989.
Central limit theorems for infinite urn models, Ann. Probab. 17,
1255--1263.

{\sc Eliazar, I.}, 2011. The Growth Statistics of Zipfian Ensembles:
Beyond Heaps' Law, Physica (Amsterdam) 390, 3189.

{\sc Gerlach, M., and  Altmann, E. G.}, 2013.
Stochastic Model for the Vocabulary Growth in Natural Languages. 
Physical Review X 3, 021006.

{\sc Gnedin, A., Hansen, B., Pitman, J.}, 2007.
Notes on the occupancy problem with
infinitely many boxes: general
asymptotics and power laws.
Probability Surveys,
Vol. 4, 146--171.

{\sc Heaps, H. S.}, 1978.
Information Retrieval: Computational and Theoretical Aspects, Academic Press.

{\sc Herdan, G.}, 1960.
 Type-token mathematics, The Hague: Mouton.

{\sc Hwang, H.-K., Janson, S.}, 2008. Local Limit Theorems for Finite and Infinite 
Urn Models. The Annals of Probability, Vol. 36, No. 3,  992--1022.

{\sc van Leijenhorst, D. C., van der Weide, T. P.}, 2005.  A Formal
Derivation of Heaps' Law, Information Sciences (NY)
170, 263.

{\sc Mandelbrot, B.}, 1965. 
Information Theory and Psycholinguistics. In: B.B. Wolman and E. Nagel. Scientific psychology. Basic Books.

{\sc Nicholls, P. T.}, 1987.  
Estimation of Zipf parameters. J. Am. Soc. Inf. Sci., V. 38, 443--445.

{\sc Petersen, A. M., Tenenbaum,  J. N.,  Havlin, S., Stanley,  H. E., Perc,  M.}, 2012.
Languages cool as they expand: Allometric scaling and the decreasing need for new words.
Scientific Reports 2, Article No. 943.

{\sc Smirnov, N.V.}, 1937. On the omega-squared distribution, Mat. Sb.  2,  973--993 (in Russian).

{\sc Zakrevskaya, N., Kovalevskii, A.}, 2019. 
An omega-square statistics for analysis of correspondence of small texts to the
Zipf---Mandelbrot law. In: Applied methods of statistical analysis. Statistical computation and simulation.
Proceedings of the International
Workshop, Novosibirsk, NSTU, 488--494.

{\sc Zipf, G. K.}, 1936. The Psycho-Biology of Language, Routledge, London, 1936.

\end{document}